\documentclass[12pt]{article}
\usepackage{amsmath}
\usepackage{amsfonts}
\usepackage{amssymb}
\usepackage{graphicx}

\newcommand{\ls}[1]
   {\dimen0=\fontdimen6\the\font \lineskip=#1\dimen0
\advance\lineskip.5\fontdimen5\the\font \advance\lineskip-\dimen0
\lineskiplimit=.9\lineskip \baselineskip=\lineskip
\advance\baselineskip\dimen0 \normallineskip\lineskip
\normallineskiplimit\lineskiplimit \normalbaselineskip\baselineskip
\ignorespaces }

\numberwithin{equation}{section}
\newtheorem{lemma}{Lemma}[section]
\newtheorem{proposition}{Proposition}[section]

\newcommand{\DD}{\mathbb{D}}%
\newcommand{\reals}{\mathbb{R}}%
\newcommand{\calB}{{\mathcal B}}
\newcommand{\calF}{{\mathcal F}}
\newcommand{\calL}{{\mathcal L}}
\newcommand{\calP}{{\mathcal P}}
\renewcommand{\tilde}{\widetilde}
\renewcommand{\hat}{\widehat}

\newcommand{\bfcdot}{{\boldsymbol \cdot}}
\newcommand{\won}{{\boldsymbol 1}}
\newcommand{\eps}{\varepsilon}

\newcommand{\pp}{\partial}
\newcommand{\om}{w}

\newcommand{\oo}{\overline}

\newcommand{\Hatx}{\,{\overline{\overline{x}}}\,}

\newcommand{\bcdot}{{\boldsymbol\cdot}}
\newcommand{\subH}{_{_{\scriptstyle{H}}}}
\newcommand{\suptwo}{^{^{\scriptstyle{2}}}}

\newcommand{\half}{\frac{1}{2}}
\newcommand{\dom}{\operatorname{dom}}
\newcommand{\trace}{\operatorname{trace}}
\renewcommand{\div}{\operatorname{div}}

\newcommand{\proof}{\noindent {\textbf{Proof:}\ }}

\newcommand{\hsp}{{\hspace*{\parindent}}}

\def\squarebox#1{\hbox to #1{\hfill\vbox to #1{\vfill}}}
\newcommand{\qed}{\hspace*{\fill}
           \vbox{\hrule\hbox{\vrule\squarebox{.667em}\vrule}\hrule}\smallskip}

\newcommand{\Pair}[4]{\raisebox{-1.2ex}{\:\mbox{\tiny$#1$}}\!
       \langle#2,#3\rangle\!\raisebox{-1.2ex}{\mbox{\tiny $#4$}\:}}
 \newcommand{\clf}[2]{\Pair{W}{#1}{#2}{W^*}}
 \newcommand{\clff}[2]{\Pair{W^*}{#1}{#2}{W}}
 \newcommand{\pair}[2]{\Pair{}{#1}{#2}{}}

\begin{document}

\title{On mutual information, likelihood-ratios and estimation error for the
additive Gaussian channel}

\author{Moshe Zakai\thanks{Department of Electrical Engineering, Technion--Israel
Institute of Technology, 
Haifa 32000, Israel.  Email: zakai@ee.technion.ac.il}}

\date{May 15, 2005} 

\maketitle

\thispagestyle{empty}

\vspace{1.3cm}
\ls{1.2}
\begin{abstract}
This paper considers the model of an arbitrarily distributed signal $x$ observed
through an added independent white Gaussian noise $w , \quad  y=x+w$.  New relations between the minimal mean
square error of the non-causal estimator and the likelihood ratio between $y$ and
$\om$ are derived.
This is followed by an extended version of a recently derived
relation between the mutual information $I(x;y)$ and the minimal mean square
error. These results are applied to derive infinite dimensional versions of
the Fisher information and the de Bruijn identity. 
A comparison between the causal and non-causal estimation errors yield a
restricted form of the logarithmic Sobolev inequality.
The derivation of the results is based on the Malliavin calculus.
\end{abstract}

\vspace{2cm}


\newpage
\ls{1.2}
\section{Introduction}

\hsp
Let $w_t , \; 0\le t \le T $ denote the standard $d$-dimensional Wiener process and $w'_t$
the related white noise.  The white noise channel is, roughly speaking,
defined by
$ y'(t)= x'(t)+w'_t$, $0\le t \le T$ where 
$x'(t)$ is a signal process independent of
the white noise process $ w'_t$.
In the context of detection theory, the key entity is $\ell(y)$, the
likelihood ratio, i.e.\ the Radon-Nikodym derivative of the measure induced by
the $\{y'(t)$,
$t\in [0,T]\}$ process with respect to the measure induced by the white noise 
$\{w'_t, t\in [0,T]\}$.  In the context of filtering theory the key entities 
are the causal
and the non causal estimates, i.e.\ the conditional mean
$E( x'_t|y'_\eta,\; 0 \le \eta \le t)$ or
$E( x'_t|y'_\eta, \; 0 \le \eta \le T)$ respectively.  
In addition to this pair
of random entities, there are also averaged entities such as the averaged minimal errors which amount to
$\int_0^T E( x'_t-E(  x'_t|y'_\eta,\; \eta\in [0,t]))^2 dt$, 
or $\int_0^T E( x'_t-E( x'_t|y'_\eta,\; \eta \in [0,T]))^2 dt$; 
and on the other
hand, the mutual information
between the paths
$\{y=(y_\eta,\; \eta \in [0,T])\}$ and
$\{x= (x_\eta,\; \eta \in[0,T])\}$, i.e.
$$
E\log \frac{dP(x;y)}{d(P(x) P(y))}
$$
where the expectation is w.r.\ to the $P(x;y)$ measure, and also the relative
entropy $E\ell(y)$.
Relations between the likelihood ratio 
$\ell(y)$ and the causal conditional expectation  
were discovered in the late 60's and this was soon followed
by a relation between the mutual information and the causal mean square error
\cite{5}, \cite{1}, \cite{beta}.  These relations which 
involved causal mean square errors were based on the Ito calculus.
Similar problems for the non causal estimator
were also considered \cite{2}, \cite{4}.
The formulation and results in the non causal case 
were restricted to the finite dimensional time discrete model of the
Gaussian channel.  Recently, however, Guo, Shamai and Verd\'u (GSV) \cite{3} 
applied information theoretic arguments to derive new 
interesting results relating
the mutual information with non causal estimation in Gaussian channels.

The Ito calculus which has proved to be a powerful tool for the relations
associated with causal estimation could not be applied to problems related to
non causal problems which explains the slow progress in the direction of
relations for
non causal estimates.  However,
the development of the Malliavin calculus, namely, the stochastic calculus of
variation which was introduced in the mid 70's led in the early 80's to
results which prove to be a very useful tool
for the non causal type of problems.

The purpose of this paper is to apply the Malliavin calculus in order to
derive the extension of the finite 
discrete time relations between non causal
estimation and likelihood ratios to continuous time (section~4), 
and to prove an extended version of the results of
\cite{3} relating the mutual information with causal estimation error
(section~5).  The modelling of the additive Gaussian channel on the abstract
Wiener space (in contrast to the $d$-dimensional Wiener process on the time
interval $[0,T]$) yield in sections~4 and 5 results of wide applicability, e.g.\
for the filtering and transmission of images and random fields. 
The relation of these results to the de Bruijn identity, causal filtering and
the logarithmic Sobolev inequality are discussed in sections~6 and 7.

In the next section we define the Abstract Wiener Space which generalizes the
classical $d$-dimensional Wiener process, and formulate the additive
Gaussian channel which will be considered in the paper. Also, the
problems considered in sections~4 and 5 are outlined in this section.
Section~3 is a very short introduction to the Malliavin calculus.  
Section~4 presents the
results relating likelihood ratios (R-N derivatives) with non-causal 
least square
estimates cf.\ remark~2 in section~4 for possible applications of these
results to nonlinear filering.  In section~5 we derive an extended version of the GSV results. 
These results are applied in section~6 to consider the notions of Fisher
information and the de Bruijn identities in an infinite dimensional setup.
Section~7 deals with abstract Wiener spaces endowed with a time parameter.   
This enables the comparison of results for causal estimations with
corresponding results for non-causal estimation.
 It is shown that a restricted form of
the logarithmic Sobolev inequality follows directly from the results derived
in this paper.

\paragraph{Acknowledgement:}
We wish to express our thanks to Shlomo Shamai for calling our attention to the problems
considered in this paper and providing us with a preliminary version of
\cite{3}, and to Suleyman \"Ust\"unel and Ofer Zeitouni for useful comments.

\section{The underlying Wiener space and the additive channel model}

\noindent
\textbf{A.}\quad Consider, first, a standard one-dimensional Wiener process on
$[0,1]$, say $\om(t), t\in [0,1]$.  Let 
$\{\eta_i (t), i=1,2, \cdots, t\in [0,1] \}$ be a complete orthonormal base on $[0,1]$.
Set $e_i(t) = \int_0^t \eta_i(s) ds$, then
$\sum_1^n \int_0^1 \eta_i(s) dw_s \cdot e_i(t)$ converges to $w(t)$
in quadratic mean.  We will denote the sequence of independent Gaussian,
identically distributed (i.i.d.) random variables 
$\left\{ \int_0^t \eta_i (s) dw(s) = \int_0^t \frac{d e_i(s)}{dt} dw(s)
\quad i=1,2,\cdots\right\}$ by $\delta e_i, \quad i=1,2,\cdots $.  Then
\begin{equation}
\label{new2.1}
w(t) = \sum_1^\infty \delta (e_i)\: e_i(t)\,.
\end{equation}
Now $\{e_i(t), t\in [0,1], i=1,2,\cdots\}$ can be considered as a 
C.O.N. base of an Hilbert
space $H$ of functions $h(t), t \in [0,1]$ with scalar product
$(h_1, h_2)_H = \int_0^1 \frac{dh_1(s)}{ds} \frac{dh_2(s)}{ds} ds$.
This space $H$ is known as the Cameron Martin space.  Note that the Wiener
process which is continuous but not differentiable is not an element in $H$.
The same notation goes over to the case of the
$d$-dimensional Wiener process with
$w(t), \eta(t), h(t), e(t)$ taking values in 
$\reals^d$ and 
\begin{equation}
\label{new2.2}
\delta e = \int_0^t \sum_1^d \frac{d}{dt} e_j (s) d w_j(s)
\end{equation}
\begin{equation}
\label{new2.3}
E \delta e =0, \quad E(\delta e)^2 = \sum_1^d \int_0^1 \left(\frac{d}{dt} e_j
(s)\right)^2 \delta s \,.
\end{equation}
In this model we will consider $y(t) = x(t) +w(t)$ where $((t), t\in [0,1])$
takes values in the Cameron-Martin Space $H$ and  $w(t)$, the Wiener process
takes values in the space of $R^d$ valued continuous function considered as a
Banach space $W$ under the norm 
$|w(t), t\in [0,1]|_W = \sup_{t\in [0,1]} |w(t)|_{\reals^d}$.

In addition to the Banach space $W$ we have to consider the space $W^*$ of all
continuous functionals on $W$ and it can be shown that $W^*$ is a dense
subspace of $H$ (cf.\ e.g.\ \cite{15}).  Hence for $e\in W^*$, it also holds
that $e\in H$ and ``$e\in W^*$ operating on $\om\in W$'' is the stochastic
integral
\begin{equation}
\label{new2.4}
\clf{w}{e} = \delta e\,.
\end{equation}
An abstract model for the Wiener process in terms of the spaces $W, W^*, H$
and the Wiener measure $\mu_W$ is considered in the next subsection.  The
reader can skip this step by interpreting the triplet $(W, H, \mu_W)$ as the
d-dimensional
Wiener process as in equations \eqref{new2.1}, \eqref{new2.2}, \eqref{new2.3}.

\bigskip
\noindent
\textbf{B.}\quad The Abstract Wiener Space (AWS) is an abstraction of this model where
$W$ is a separable Banach space and $H$, the Cameron-Martin space, is an
Hilbert space densely and continuously embedded in $W$.  The dual space to $W$
(the space of continuous linear functionals on $W$) is denoted $W^*$ and
assumed to be continuously and densely embedded in $H$.  The Abstract Wiener
Space $(W,H,\mu_W)$ supports a $W$-valued random variable $w$ such that for every
$e\in W^*$, $\:\delta e:= \clf{w}{e}$ is a $N(0, |e|\subH^2)$ random
variable.
Cf.\ e.g.\ \cite{11}, or appendix~B of \cite{12} 
and the references therein,
for further information on the AWS.
Note that, unlike the classical case, the Abstract Wiener Space  
does not have any time-like parameter
(this however can be added cf.\  section~7).

\bigskip
\noindent
\textbf{C.}\quad In order to introduce the general setup of the additive 
Gaussian channel, let 
$(W,H,\mu_W)$ be an abstract Wiener space and let $(H, \sigma (H), \mu_X)$ 
be
a probability space on the Cameron-Martin Hilbert space $H$ which is induced by an
$H$-valued r.v.\ $X$.
Let $\theta = (x,w), x\in H$ and $\om\in W$, set $\Theta =\{\theta\}$ and
consider the combined probability space
\begin{equation}
\label{2.3}
(\Theta, \calF, \calP) = \Bigl(\Theta, \sigma (H) \vee \sigma (W), \mu_X
\times \mu_W\Bigr)
\end{equation}
which is the space of the mutually independent `signals' $x$ and `noise' $w$.
Now, since $H$  
is continuously embedded in $W$ we can identify $x$ with its image
in $W$ and defined the additive Gaussian channel as
\begin{equation}
\label{2.4}
y(\theta) = \rho x + w \, ,
\end{equation}
where $\rho$ is a free scalar `signal to noise' parameter which will become
relevant in Section~5.  We will denote by $X$ and $Y$ the sigma fields
induced on $W$ by the r.v.'s $x$ and $y$ respectively.
Note that $y$ and $w$ are $W$ valued, $x$ is $H$ valued  and
we identify $x$ with its image in $W$.
In fact we will make throughout this paper, just for reasons of simplicity, the additional assumption that
$x$ is $W^*$ valued. As mentioned earlier, since 
$W^* \subset H \subset W$ we can also consider $x$ to be $H$ or $W$ valued.

In section~4 we will be interested in the relation between two types of
objects.  The first class of objects is
$$
E\Bigl( (x,e)\subH\:|Y\Bigr) \quad \text{\ and \ } \quad 
E\Bigl((x,e_1)\subH \cdot (x,e_2)\subH \: | Y\Bigr)
$$
for $e, e_1, e_2\in W^*$ or globally 
\begin{equation}
\label{2.5}
\Hatx =E(x|Y) \quad \text{\ and \ } \quad  
\overline{\overline{(x,x)}}\! \subH = E\Bigl((x,x)\subH|Y\Bigr)\,.
\end{equation}
The second class of objects are the likelihood ratio (the R-N derivative) between the
measures induced by $y$ and the one induced by $w$ on $W$.  This likelihood
ratio will be denoted $\ell(\om)$, $\om\in W$.  Note that 
if $W$ is infinite dimensional then
the measure induced
by $x$ is singular with respect to the measure induced by $w$, (since $x\in H$
while $w\not\in H$).

In section~5 we will consider the relation between $I(X;Y)$ or rather
$dI(X;Y)/d\rho$ and the non-causal filtering error:
$$
\Bigl\{
E|x|\subH\suptwo - \Bigl|E(x|Y) \Bigr|\subH\suptwo \Bigr\} = 
E\Bigl\{|x|\subH\suptwo - |\Hatx |\subH\suptwo \Bigr\}
\,.
$$
A related result for $d(E\log \ell(\om))/d\rho$ is considered in section~6 and shown to be an extended
version of the De Bruijn identity.

\section{A short introduction to the Malliavin calculus}

\hsp
For further information cf.\ e.g.\ \cite{6}, \cite{8}, \cite{11} or appendix~B
of \cite{12}.

\noindent
\textsl{(a) The gradient}

Let $(W,H,\mu)$ be an AWS and let $e_i, i=1,2,\dots$ be a sequence of elements
in $W^*$.  Assume that the image of $e_i$ in $H$  
form a complete orthonormal base in $H$.  Let $f(x_1, \dotsc, x_n)$ be a
smooth function on $\reals^n$ and denote by $f'_i$ the partial derivative of
$f$ with respect to the $i$-th coordinate and let $\delta e$ be as discussed
in the previous section.

For  cylindrical
smooth random variables $F(\om) = f(\delta e_1, \dotsc, \delta e_n)$, define\\
$\nabla_h F = \left.\frac{dF(w+\eps h)}{d\eps}\right|_{\eps=0}$.  Therefore we
set the following: $\nabla_hF=(\nabla F, h)$ where
$\nabla F $, the gradient, is $H$-valued.  For $ F(w) = \delta e , \, \nabla F = e $, and
\begin{equation}
\label{3.1}
\nabla F = \sum_{i=1}^n f_i' (\delta e_1,
\dotsc, \delta e_n) \cdot e_i
\,.
\end{equation}
It can be shown that this definition is closable in $L^p(\mu)$ for any $p>1$,
which means that it can be extended to a wider class of functional as we will
see below.
We will restrict ourselves to $p=2$,
consequently the domain of the $\nabla$ operation can be extended to all
functions $F(\om)$ for which
there exists a sequence of smooth cylindrical functions
$F_m$ such that $F_m\to F$ in $L_2$ and $\nabla F_m$ is Cauchy in 
$L_2(\mu, H)$.  In this case set $\nabla F$ to be the $L_2(\mu,H)$ limit of
$\nabla F_m$.  This class of r.v. will be denoted $\DD_{2,1}$. It is a closed
linear space under the norm
\begin{equation}
\label{3,2}
\|F\|_{2,1} = E_\mu^{\half} (F)^2 + E_\mu^{\half} |\nabla F|\subH\suptwo
\,.
\end{equation}
Similarly let $K$ be an Hilbert space and $k_1,k_2, \dotsc$ a complete
orthonormal base in $K$.  Let $\varphi$ be the smooth $K$-valued function
$ \varphi =\sum_{j=1}^m f_j(\delta e_1, \dotsc, \delta e_n) k_j$ define
\begin{equation}
\label{3.3}
\nabla \varphi = \sum_{j=1}^m 
\sum_{i=1}^n (f_j)'_{i} (\delta e_1, \dotsc, \delta e_n) e_i
\otimes k_j
\end{equation}
and denote by $\DD_{2,1}(K)$ the completion of $\nabla \varphi$ under the norm
\begin{equation}
\label{3.4}
\|\varphi \|_{2,1} = E^{\half} \Bigl( |\varphi|_K^2 +  | \nabla \varphi
|_{H\times K}^2 \Bigr)\,.
\end{equation}
Note that this 
enables us to define recursively $\nabla^nF(\om)$ for $n>1$.

\noindent
\textsl{(b) The divergence (the Skorohod integral)}

A few introductory remarks.  Let $v(x), x\in \reals_n$ take values in
$\reals_n$, $v(x) = \sum_1 v_i(x) \rho_i$, where the $\rho_i$ are orthonormal
vectors in $\reals_n$.  Assume that the $v_i$ and $F(x)$ are smooth and
converge ``quickly enough'' to zero as $|x|\to\infty$.  Then the following
``integration by parts formula'' holds
\begin{equation}
\label{3.4+1}
\int_{\reals_n} (v(x), \nabla F(x)) dx = - \int_{\reals_n} F(x) \div v \,dx\,,
\end{equation} 
where div is the divergence:
$$
\div v = \sum_1^n \frac{\pp v_i}{\pp x_i}\,.
$$
Note that the gradient and divergence are differential operations,
and
equation~\eqref{3.4+1} deals with integration with respect to the Lebesgue
measure on $\reals_n$.
In this subsection we are looking for an analog of the divergence operation on
$\reals_n$ which will yield an integration by parts formula with respect to
the Wiener measure.

Let $u(\om)$ be an $H$-valued r.v.\ in $(W,H,\mu)$, $u$ will be said to be in  
$\dom_2 \delta $ if $E|u(\om)|\subH\suptwo < \infty$ and there exists a r.v.\ say 
$\delta u$ such that for all smooth functionals 
$f(\delta e_1, \dotsc, \delta e_n)$ and 
all $n$
the ``integration by parts'' relation
\begin{equation}
\label{3.5}
E\Bigr(\nabla f, u(\om)\Bigr)\subH = E(f \cdot \delta u)
\end{equation}
is satisfied.
$\delta u$ is called the divergence or Skorohod integral. 
A necessary and sufficient condition for a square integrable $u(w)$ to be in
$\dom_2\delta$ is that for some $\gamma=\gamma(u)$,
$$
\Bigl|E(u(w), \nabla f)\subH\Bigr| \le \gamma E^{\half} f^2 (w)
$$
for all smooth $f$.
 Note that while the
definition of $\nabla f$ (at least for smooth functionals) is invariant under
an absolutely continuous change of measure, this is not the case for the
divergence which involves expectation in the definition.  For non-random 
$h\in W^*$, $\delta h = \pair{h}{w}$, setting
$f=1$ in \eqref{3.5} yields that $E\delta h=0$.  It can be shown that if
$u\in \DD_{2.1} (H)$ then $u\in\dom_2 \delta$.
Also, for smooth $f(\om)$ it can be
verified directly that
$$
\delta \Bigl( f(\om) h\Bigr) = f(\om) \delta h - (\nabla f, h)\subH
$$
and more generally under proper restrictions
\begin{equation}
\label{3.6}
\delta \Bigl(f(\om) u (\om)\Bigr)
= f(\om) \delta u - \Bigl(  \nabla f, u(\om)\Bigr)_H
\,.
\end{equation}
Consequently, if $E|u|\subH\suptwo < \infty$,  and 
$ \nabla u $ is of $\trace$ class then
\begin{equation}
\label{new3.7}
\clff{u}{\om} = \delta u + \trace \nabla u \, .
\end{equation}
where
for an operator $A$ on $H$ and 
$e_i, i=1,2,\dots$ a CONB on $H$, define 
$$
\trace A = \sum_1^\infty (e_i, Ae_i)
$$
provided the series converges absolutely and in this case $A$ is said to be of
trace class.  
Among the interesting facts about the divergence operator, let us also note that for
the classical Brownian motion and if
$$
(u(\om))(\cdot) = \int_0^\bcdot u'_s(\om) ds
$$
and $u'_s(\om)$ is adapted and square integrable then
$\delta u$ coincides with the Ito integral i.e.\
$\delta u=\int_0^1 u'_s(\om) dw_s$.

\medskip
\noindent
\textsl{(c)}\quad Let $(W,H,\mu)$ be an abstract Wiener space and let $\mu_1$
be another probability measure on the same space $(W,\sigma \{W\})$.
Assume that $\mu_1$ is absolutely continuous with respect to $\mu$.  Set
$$
\ell(\om) = \frac{d\mu_1}{d\mu} (\om) \quad \text{\ and \ } \quad
Q(\om) = \{\om: \ell (\om) > 0\}
$$
$E_1$ and $E_0$ will be used to denote the expectation with respect to the
measures $\mu_1$
and $\mu$ respectively.  We will use the convention $0\log 0=0$ throughout the
paper.

Following the definition in 3(b), we define the divergence with respect to
$\mu_1$ to be as follows.  The $H$-valued random variable $u(\om)$ will be
said to be in $\dom_2^1 \tilde \delta$ if there exists a r.v., 
say $\tilde\delta u$, which is $ L^2 $ under $ u_1 $ and such that for all smooth  r.v.s
$f(\om)$, it holds that
$$
E_1\Bigl( f(\om) \cdot \tilde \delta u\Bigr)= E_1 (u, \nabla f)\subH\,.
$$
The relation between $\tilde \delta u$ and $\delta u$ is given by the
following lemma.
\begin{lemma}
\label{lem-31}
Assume that $ \ell(\om) \in \DD_{2,1}$, $u \in \dom_2 \delta $,
$\ell \cdot \delta u \in L_2$ and $\ell \cdot \nabla u \in \DD_{2,0} (H)$ and
$\mu_1 \ll \mu_W$ 
(where $\DD_{2,0}(H)$ is the completion of \eqref{3.3} under the $H$-norm).
Then $u  \in \dom_2^1 \tilde \delta$ and
\begin{equation}
\label{3.7}
\tilde\delta u = \won_Q(\om) (\delta u - \Bigl(\nabla \log \ell(\om), u
(\om)\Bigr)\subH
\end{equation}
\end{lemma}

\proof
Since $f(\om)$ is a  smooth r.v., 
$\ell\cdot f \cdot \delta u - f(\nabla f, u)\subH$ is in $L_1$ and
$\ell(\om) \nabla \log \ell(\om) = \nabla \ell(\om)$ a.s.-$\mu$.  Hence
\begin{align*}
E_1 \Bigl(f(\om) \delta u - f(\om) \Bigl( \nabla \log \ell (\om), u)\subH \Bigr)
&= E_0\Bigl(\ell \cdot f\cdot \delta u - \ell f (\nabla \log \ell, u)\subH\Bigr)\\
& = E_0\Bigl\{(\nabla(\ell\cdot f), u\Bigr)\subH - f(\nabla \ell, u) \Bigr\} \\
& = E_0 \Bigl( \ell(\nabla f, u)\subH\Bigr)\\
& = E_1 (\nabla f, u)\subH \tag*{\qed} 
\end{align*}

\section{Relations between the estimation error and the likelihood ratio}

\hsp
Let $(W,H,\mu), (H, \sigma (H), \mu_X)$, $(\Theta, \calF, \calP)$ and
$y(\theta) = \rho x + w$ be as in section~2.  We will further assume that the
$H$-valued r.v. $x$ is actually $W^*$ valued, and 
$\exp \alpha (x,h)\subH \in L^1 (\mu_X)$ for all real 
$\alpha $ and all $h\in W^*$. 
The measures induced by $y$ and
$x$ on $W$ will be denoted $\mu_Y$ and $\mu_X$
respectively.  The conditional probability 
induced on $W$ by $y(\theta)$ conditioned on $x$ will be denoted by 
$\mu_{Y|X}$.  Similarly, $\mu_{X|Y}$ will denote the conditional probability
induced on $W^*$ of $x$ conditioned on $y$ (cf.\ e.g.\ \cite{14} for the
existence of these conditional probabilities).

By the Cameron-Martin theorem (cf.\ e.g.\ \cite{12}) and since $x$ and $w$
are independent, we have
\begin{equation}
\label{4.2}
\frac{d\mu_{Y|X}}{d\mu_W} (\om) = \exp
\left( \rho \pair{w}{x} - \frac{\rho^2}{2} |x|\subH\suptwo\right), \quad w \in W
\end{equation}
which by our assumptions belongs to $L_p$ for all $p>0$.  Hence, denoting by
$\mu_X(dx)$ the restriction of $\calP$ to $H$:
\begin{align}
\label{4.3}
\ell(\om) & = \frac{d\mu_Y}{d\mu_W}(\om) =\int_H\, \frac{d\mu_{Y|X}}{d \mu_W}
(\om,x)
\mu_X(dx) \notag\\
& = \int_H \exp \left(\rho \pair{w}{x} - \frac{\rho^2}{2} |x|\subH\suptwo\right)
\mu_X(dx)
\end{align}
\begin{proposition}
\label{prop-41}
Under these assumptions it holds that\\
(a)

\begin{equation}
\label{4.4}
\begin{aligned}
(\nabla \ell, h)\subH =:
&\nabla_h \ell(w) = \rho \ell(w) (\Hatx,h)\subH, \quad \forall h \in H \quad
\mathrm{hence:}  \\
&\nabla \ell = \rho \ell(w) \Hatx
\quad \text{or} \quad 
\Hatx  = \frac{1}{\rho} \nabla \log \ell(\om)
\end{aligned}
\end{equation}
a.s.\ $\mu_W$.  Note that $\oo{\oo a}$
denotes the conditional expectation  conditioned on $Y$, equation~\eqref{2.5}.
\\
(b)

\begin{equation}
\label{4.5}
\Bigl(\nabla^n \ell(\om), h_1 \otimes \dots \otimes h_n\Bigr) 
_{H^{\otimes n}}
= \ell (\om) \rho^n \overline{\overline{\left(\prod_{i=1}^n (h_i, x)\right)}}
\end{equation}

\noindent 
(c)~in particular $ \trace \nabla^2 \ell (w) $ exists and a.s.\ $\mu_W$
\begin{equation}
\label{4.6}
\nabla_{h_1,h_2}^2 \ell(\om) = \rho^2 \ell(\om)\; \overline{\overline{\Bigl((h_1,x) \cdot 
(h_2, x)\Bigr)}}
\end{equation}
and
\begin{equation}
\label{4.5b}
\nabla_{h,h}^2 \log \ell(\om) = \rho^2 
\Bigl(\;\overline{\overline{\Bigl((x,h)^2\Bigr)}} - (\Hatx ,h)^2\Bigr)
\end{equation}
where $\nabla_{h_1, h_2}^2 \varphi =: (\nabla (\nabla \varphi, h_2),
h_1)\subH$, cf.\ also \eqref{4.8}

\noindent
(d)

\begin{equation}
\label{4.7}
\oo{\oo{\Bigl(\prod_{i=1}^n(h_i,x)\Bigr)}} =
(h_n, \Hatx)\; \oo{\oo{\Bigl(\prod_1^{n-1} (h_i, x)\Bigr)}} +\nabla_{h_n}\;
\oo{\oo{\Bigl(\prod_{i=1}^{n-1} (h_i, x)\Bigr)}}\, .
\end{equation}
\end{proposition}

\paragraph{Remark 1:}
Let $E_1$
denote the measure induced by $y$ on $W$ and let $E$ denote expectation w.r. to the measure in \eqref{2.3}.
For an operator $A$ on $H$ and 
$e_i, i=1,2,\dots$ a CONB on $H$, define 
$$
\trace A = \sum_1^\infty (e_i, Ae_i)
$$
provided the series converges.  
Consequently, we have from \eqref{4.5b} and \eqref{4.4}
that 
\begin{equation}
\label{4.8}
\begin{split}
E_1 \trace \nabla^2 \log \ell (w) = \rho^2 E |x-\Hatx|\subH\suptwo
&= \rho^2 \Bigl(E|x|\subH\suptwo - E|\Hatx|\subH\suptwo\Bigr) \\
& = \rho^2 E|x|\subH\suptwo - E|
\nabla \log \ell (w) |\subH\suptwo\,.
\end{split}
\end{equation}
(c.f.\ also equations \eqref{6.3} and \eqref{6.4}).

\paragraph{Remark 2:}
(a)~Consider the case where the abstract Wiener space is a classical 
Wiener space $\reals^n$, then $u\in H$ is of the form
$\int_0^\cdot u' (s) ds$, $x\in H$ is of the form
$\int_0^\cdot x' (s) ds$ where $x'(s) \in \reals_n$ and
$\int_0^T |x'(s)|_{R_n}^2 ds < \infty$.
Further assume that $E\int_0^T |x'(s)|_{\reals_n}^2 ds < \infty$ and
$x'(s)$ is a.s.\ continuous on $[0,T]$.  Then given some $t\in [0,T]$, one can
consider a sequence of linear functionals $h_n$ such that $(h_n,x)$ converges
in  $L^2$ to $x(t,\om)$ and extend the results of proposition~\ref{prop-41} to
$\Hatx_t := E(x'(t)|Y)$ for any $ t\in [0,T] $.  
(b)~Consider equation \eqref{4.3}, given $\om \in W$ we can replace the
integration with respect to $\mu_X(dx) $ with a Monte-Carlo approximation.
Similarly we can replace $\langle \nabla^n \ell (\om), h_1 \otimes \dots \otimes
h_n\rangle$ by applying $\nabla^n$ to the integrand of equation~\eqref{4.3} and
then replacing again the integrand with a Monte-Carlo approximation.  This can
then be applied to derive a numerical approximation to 
$\overline{\overline{\prod_{i-1}^n (h_i, x)}}$ i.e.\ the non-adapted non-linear fitlering of 
$\prod_{i=1}^n (h_i, x)$.
We will not follow these directions.

The following lemma will be needed in the proof of Proposition~\ref{prop-41}:
\begin{lemma}
\label{lem-41}
Assume that $\mu_Y$ and $\mu_{Y|X}$ are absolutely continuous with respect to 
$\mu_W$ then for all bounded and measurable functions $\psi$ on $\Theta$
$$
\int_{X\times W}
\psi(x,y) \frac{d\mu_{Y|X}}{d\mu_W}(y,x) \mu_X (dx) \times
\mu_W (dy) =
\int_{X\times W} \psi(x,y)  \frac{d\mu_Y}{d\mu_W}(y)
\mu_{X|Y} (dx;y) \mu_W (dy)
\,.
$$
\end{lemma}
\textbf{Proof of Lemma:}
Let
$$
L= \int_{X\times Y} \psi (x,y) \mu_{X,Y} (dx,dy)
\,.
$$
Then, by Fubini's theorem
\begin{align*}
L & = \int_{X\times Y} \psi(x,\om) \mu_{Y|X} (d\om,x) \mu_X(dx)\\
& = \int_{X\times Y} \psi(x,y) \frac{d_{Y|X}}{d\mu_W}(y,x)
\mu_X (dx) \mu_W (dy) \, .
\end{align*}
Since the conditional probability $ \mu_{X|Y} $ is regular (cf e.g. theorem~10.2.2 of \cite{14}) we also have
\begin{align*}
L & = \int_{X\times Y} \psi(x,y) \mu_{X|Y} (dx,y) \mu_Y (dy) \\
& = \int_{X\times Y} \psi(x,y) \frac{d\mu_Y}{d\mu_W} (y)
\mu_{X|Y} (dx,y) \mu_W (dy)\,. \tag*{\qed}
\end{align*}
\textbf{Proof of Proposition:}
From \eqref{4.2} and \eqref{4.3} and since by our assumptions we may (by
dominated convergence) interchange the order of integration and
differentiation
\begin{align*}
\nabla_h \ell(\om) 
&= \int_H \rho (\nabla_h \pair{\om}{x}) \exp\left(\rho \pair{\om}{x}-
\frac{\rho^2}{2} |x|\subH^2\right) \mu_X (dx)\\
&
= \int_H \rho(h,x) \exp \left(\rho\pair{\om}{x} - \frac{\rho^2}{2}
|x|\subH\suptwo\right) \mu_X(dx)\\
& = \int_H \rho (h,x) \frac{d\mu_{Y|X}}{d\mu_W}(\om ,x) \mu_X(dx)
\,.
\intertext{Thus, by Lemma~\ref{lem-41}}
\nabla_h \ell (w) & = \int_X \rho(h,x) \frac{d\mu_{Y}}{d\mu_W}(\om)\mu_{X|Y} (dx,w)\\
& = \rho\ell(\om) (h, \Hatx)
\end{align*}
proving \eqref{4.4}.  The same arguments also hold for repeated
differentiation
\begin{align*}
\Bigl( \nabla^n \ell(\om), h_1 \otimes \dots \otimes h_n\Bigr)_{H^{\otimes n}} 
& = \rho^n \int_X (h_1,x) \dots (h_n, x) \frac{d\mu_{Y|X}}{d\mu_W}(\om ,x)
\mu_X(dx)\\
& = \rho^n \int_x \left( \prod_{i=1}^n (h_i, x) \right)
\frac{d\mu_{Y}}{d\mu_W}(\om) \mu_{X|Y} (dx,w) \, ,
\end{align*}
which yields \eqref{4.5}.  \eqref{4.6} follows directly from \eqref{4.5}
since
\begin{equation}
\label{4.9}
\nabla_{h_1,h_2}^2 \ell (\om) = \rho^2 \ell(\om) 
\;\oo{\oo{\Bigl((h_1,x) \cdot (h_2 \cdot x)\Bigr)}}
\end{equation}
therefore
\begin{align}
\nabla_{h,h}^2 \log \ell(\om) & = \frac{1}{\ell(\om)} \nabla_{h,h}^2
\ell(\om) - \frac{1}{(\ell(\om))^2} (\nabla_h \ell(\om))^2\notag\\
& = \frac{1}{\ell(\om)} \nabla_{h,h}^2 \ell(\om) - |\nabla_h\log\ell (w)|^2
\tag{4.9a}
\end{align}
proving \eqref{4.6} and \eqref{4.5b}.  From \eqref{4.5} we have
\begin{align*}
\rho^n \ell(\om)\;
\oo{\oo{\Bigl(\prod_{1}^{n} (h_i, x)\Bigr)}}
&= \nabla_{h_n} \Bigl( \nabla^{n-1} \ell(\om), h_1 \otimes \dots \otimes
h_{n-1}\Bigr)_{H^{\otimes (n-1)}} \\
& = \nabla_{h_n} \left(\ell(\om) \cdot \rho^{n-1} \cdot
\oo{\oo{\Bigl(\prod_{1}^{n-1} (h_i, x)\Bigr)}}\;
\right) 
\end{align*}
and \eqref{4.7} follows.
\qed

We conclude this section with some results for $\delta\Hatx$ and
$\tilde\delta\Hatx$ (cf. part (c) of section~3).  By the assumptions of this
section
$\Hatx\in \dom_2\delta$ and $\Hatx \in \dom_2^1 \delta$.  Therefore by
\eqref{4.4}
$$
\delta \Hatx = \frac{1}{\rho} \delta \nabla \log \ell(\om)
\,.
$$
Note that $\calL = \delta \nabla$ is the number operator, i.e. if $ \alpha (w) $ is a square integrable r.v.
of the Wiener space and $ \alpha = \sum_{n=1} \, I_n $,
where $ I_n $ is the Wiener chaos decomposition
of $ x$; then, formally, $ \calL\alpha =\sum_{n=1} \, n I_n $.
Therefore if
$\calL\alpha(\om) \in L_2$ and $E(\calL\alpha(\om)) =0$ then 
$\calL^{-1}\calL \alpha$ is well
defined, consequently it holds by equation \eqref{4.4} that
\begin{equation}
\label{4.10}
\ell(\om) = c \cdot \exp \rho \calL^{-1} \delta \Hatx\,.
\end{equation}
where $c$ is a normalizing constant.
For $\tilde\delta \Hatx$ we have
\begin{lemma}
\label{lem-42}
$$
\tilde\delta \Hatx = \frac{1}{\rho \ell(\om)} \delta \nabla \ell(\om)
\,.
$$
\end{lemma}
\proof
By \eqref{4.4}
\begin{align*}
\delta \nabla \ell & = \rho\delta (\ell(\om) \Hatx) \\
& = \rho \ell(\om) \delta \Hatx - \rho (\Hatx, \nabla \ell)\subH\\
& = \rho \ell(w) \delta \Hatx - \rho^2 \ell(\om) (\Hatx, \Hatx)\subH
\end{align*}
and
$$
\delta \Hatx = \frac{1}{\rho\ell(\om)} \delta\nabla \ell(\om) + 
\rho (\Hatx, \Hatx)\subH
\,.
$$
Hence by Lemma~\ref{lem-31}
\begin{align*}
\tilde\delta \Hatx & = \delta \Hatx - \Bigl( \nabla \log \ell(\om), \Hatx)\subH\\
& = \delta\Hatx - \rho |\Hatx|\subH\suptwo\\
& = \frac{\calL \ell(\om)}{\rho \ell(\om)} \tag*{\qed}
\end{align*}

\section{The GSV relation between the mutual information and the mean square
of the estimation error}

\hsp
Consider the setup and assumptions in the first paragraph of section~4.
The mutual information between $x$ and $y$ is defined as
$$
I(X;Y) = \int_{X\times W} \log
\frac{d\mu_{X;Y}}{d(\mu_X \times \mu_Y)} (x,y)\; 
\mu_{X,Y} (dx,dy)\,.
$$
$E$ will denote expectation w.r. to the measure in \eqref{2.3}, (cf.\ e.g.
\cite{9}). $ E_0$ will denote expectation w.r. to the
Wiener measure and $ E_1 $ will denote expectation w.r. to the measure on $W$ induced by $y$ (hence $ Ef(y)
= E_1 f(w) = E_0\ell (w) f(w)$).
\begin{proposition}
\label{prop-51}
Under the assumptions of the previous section, it holds that
\begin{align}
\label{5.0}
\frac{dI(X;Y)}{d\rho} &= \rho E \Bigl(|x|\subH\suptwo - | \Hatx|\subH\suptwo\Bigr)\\
& = \rho E |x-\Hatx|\subH\suptwo \notag
\end{align}
\end{proposition}
\proof
By our assumptions, and since
$\frac{d\mu_{X,Y}}{d(\mu_X \times \mu_Y)} = \frac{d\mu_{Y|X}}{d\mu_W} \cdot
\frac{d\mu_W}{d \mu_Y}$, we have
\begin{align*}
I(X;Y) & = \int_{X\times W}
\left\{ \log \frac{d\mu_{Y|X}}{d\mu_W} (x,y) - \log \frac{d\mu_Y}{d\mu_W}
(y)\right\} \mu(dx,dy) \\
& = E \left(\rho\pair{y}{x} - \frac{\rho^2}{2} |x|\subH\suptwo\right) - E \log
\ell(\om)\,.
\end{align*}
Note that $E \rho\pair{y}{x} = \rho^2 E |x|\subH\suptwo$, hence
\begin{equation}
\label{5.2}
I(X;Y) = \frac{\rho^2}{2} E |x|\subH\suptwo - E_1 \log \ell(\om)
\end{equation}
and
\begin{align}
\label{5.3}
\frac{dI(X;Y)}{d\rho} & = \rho E |x|\subH\suptwo - \frac{d}{d\rho}
E_0 \ell(\om) \log \ell(\om)\\
\label{5.4}
& = \rho E |x|\subH\suptwo - E_0 \left(\frac{d\ell(\om)}{d\rho} \cdot \log
\ell(\om)\right) - 0\,.
\end{align}
Now,
\begin{align*}
\frac{d\ell(\om)}{d\rho} & = \int_X \Bigl( \pair{x}{\om} - \rho
|x|\subH\suptwo\Bigr) \frac{d\mu_{Y|X}(x)}{d\mu_W} (\om) \mu_X (dx)\,.\\
\intertext{By lemma~\ref{lem-41}}
\frac{d\ell(\om)}{d\rho} & = \int_X \Bigl( \pair{x}{\om} - \rho |x|\subH\suptwo\Bigr)
\frac{d\mu_{Y}}{d\mu_W} (\om) \mu_{X|Y} (dx) \\
& = \Bigl( \pair{\Hatx}{\om} - \rho\; \oo{\oo{(|x|\subH\suptwo)}}\; \Bigr) \ell (\om) \, .
\end{align*}
Substituting in \eqref{5.4} yields
\begin{equation}
\label{5.5}
\frac{dI}{d\rho} = \rho E|x|\subH\suptwo - E_0
\left\{\Bigl( \pair{\Hatx}{\om} - \rho\; \oo{\oo{(|x|\subH\suptwo)}}\;\Bigr)
\ell(\om) \log \ell(\om) \right\} \, .
\end{equation}
Now, by \eqref{4.4}
\begin{align*}
E_0 \ell \log \ell \pair{\Hatx}{\om}
 = & \quad E_0 \left(\frac{1}{\rho} \log \ell \pair{\nabla \ell}{\om}\right)\\
 \overset{\text{by \eqref{new3.7}}}{=} & \quad E_0\: \frac{1}{\rho}\: 
\log \ell \Bigl(\delta \nabla \ell + \trace \nabla^2
\ell\Bigr) \\
 \overset{\text{by \eqref{3.6}}}{=} & \quad E_0 \:\frac{1}{\rho}\: \delta (\log \ell \nabla \ell) -
E_0\: \frac{1}{\rho}\: \Bigl(\nabla \ell, \nabla \log \ell\Bigr)
+ E_0\: \frac{1}{\rho}\: (\log\ell\trace\nabla^2\ell) \\
 \overset{\text{by \eqref{3.5} and \eqref{4.9}}}{=} & \quad 0 - \frac{1}{\rho} E_0 \frac{1}{\ell}
(\nabla \ell, \nabla \ell) + E_0 \rho 
\left(\ell\log  \ell\;\oo{\oo{(|x|\subH\suptwo)}}\;\right)\\
 = &\quad \rho E_0 \ell (w) (\Hatx, \Hatx) + E_0 \rho \ell (\om) \log \ell (\om)
\;\oo{\oo{(|x|\subH\suptwo)}}\;\,.
\end{align*}
Substituting into \eqref{5.5} yields
$$
\frac{dI}{d\rho} = \rho E|x|\subH\suptwo -\rho E |\Hatx|\subH\suptwo + E_0 \rho  \ell
 \log \ell\;\oo{\oo{ (|x|\subH\suptwo)}}\;
- E_0 \rho \ell  \log \ell\;\oo{\oo{ (|x|\subH\suptwo)}}\;
$$
\qed

\paragraph{Remark (a):}
$E_1\log\ell(w) (=E \log\ell(y))$ is the relative entropy (or $I$-divergence or
Kullback-Leibler number) of $\mu_Y$ with respect to $\mu_W$ (cf.\ e.g. \cite{9} or
\cite{alpha}).  Equation \eqref{5.2} relates this relative entropy to the mutual
information $I(x;y)$ for the additive Gaussian channel.  By equations
\eqref{5.3} and \eqref{5.0}
\begin{equation}
\label{7.1}
\frac{d}{d \rho} E_1 \log \ell (w) = \rho E_1 | \Hatx |\subH\suptwo\,.
\end{equation}

\paragraph{Remark (b):}
Consider the following generalizations to the additive Gaussian channel.
Let $M$ be ``the space of messages which generate the signals'' $x$, i.e. $ ( M, \calB_\mu , \calP_\mu ) $ is
a probability space and $ x=g(m) \, , \quad m\in M $, where $g$ is a measurable from $ (M, \calB_\mu ) $ to
$H$. Then obviously $ I(M,Y) = I(X;Y) $. More generally, consider the case where $x$ and $m$ are related by
some joint probability on $ M\times H $ and $ w$ and $m$ are conditionally independent conditioned on $x$.
The extension of proposition~\ref{prop-51} in this context follows along the same arguments as in theorem~13 of \cite{3}
and therefore omitted.

\section{An extended version of the De Bruijn identity}

\hsp
The Fisher information matrix $J$ associated with a smooth probability density
$ p( y_1, \dotsc ,y_n ) $, $ y\in\reals_n $ is defined as  
$$
J=\left(\frac{\pp^2 \log p(y_1 ,\dotsc ,y_n)}{\pp y_i \pp y_j}\right)_{1\le i,j\le n}
$$
and then
the Fisher information which is defined by the r.h.s. of \eqref{6.1}
satisfies:
\begin{equation}
\label{6.1}
E\trace J =- E \Bigl\{\Bigl|\nabla \log p\Bigr|_{\reals^n}^2\Bigr\} \, ,
\end{equation}
where $E$ is the expectation with respect to the $p$ density.
The De Bruijn identity (cf. \cite{20} or \cite{3} and the references therein) deals with the case where $ y=x+\sqrt{t}w$
where $ w=w_1 , w_2 ,\dotsc ,w_n $ and the $ w_j \, , \quad j=1 , \dotsc ,n $ are i.i.d. $ N(0,1) $ and
$x$ is an $ \reals_n $ random variable independent of $w$. It states that
\begin{equation}
\label{6.2}
\frac{d}{dt} \, E\log p(y) = \half\, E\Bigl\{\Bigl| \nabla \log p(y)\Bigr|_{\reals^n}^2\Bigr\} \, .
\end{equation}
The Fisher information matrix cannot be extended directly to the case where $y$ is infinite dimensional.
However, the results of sections~4 and 5 yield some similar relations. Under the assumptions of section~5,
comparing \eqref{5.0} with \eqref{5.3} we have
\begin{align}
\label{6.3}
\frac{d}{d\rho}\, E_1 \log \ell(w) & = \rho E_1 |\Hatx|_H^2\notag\\
& = \frac{1}{\rho}\, E_1 |\nabla \log\ell (w)|_H^2 \, , 
\end{align}
which is ``similar'' to \eqref{6.2} and may be considered an extended De Bruijn identity.
Note that $ E_1 \log\ell (w) $ is the relative entropy of $ \mu_Y $ relative to $ \mu_W $, also, note the
difference between the $ \rho $ and the $t$ parametrizations.

Comparing \eqref{5.0} with \eqref{4.8} yields
\begin{align}
\label{6.4}
E_1 \trace \nabla^2\log \ell (w) & = \rho\, \frac{dI(x;y)}{d\rho}\notag\\
& = \rho^2 E |x|_H^2 - E_1 |\nabla \log \ell (w)|_H^2 \, ,
\end{align}
which is different from \eqref{6.1} by the $ \rho^2 E |x|_H^2 $ term.
Note that the validity of \eqref{6.4} is restricted to the case where
$\ell(w)$ is induced by a signal plus independent noise model and
\textsl{not\/} for \textsl{any\/} $\ell(w)$ which is a negative r.v.\ and whose
expectation is 1, cf.\ the concluding lines of the next section.

\section{Adding a ``time parameter'' to the abstract Wiener space}

Given an abstract Wiener space $(W, H, \mu)$ we can introduce the notion
of continuous time on this space as follows.  Let
$\{\pi_\theta, 0 \le \theta \le 1\}$ be a continuous , strictly increasing,
resolution of the identity on $H$ with
$\pi_0 = 0, \pi_1=I$.  Set 
$\calF_\theta = \sigma \{\delta \pi_\theta h, h \in H\}$ and $\calF_\bfcdot$
will denote the filtration induced by $\calF_\theta$ 
on $[0,1]$.  An $H$-valued r.v.\ $u(w)$ will be said to be adapted to
$\calF_\bfcdot$ if $(u, \pi_\theta, h)\subH $ is $\calF_\theta$ measurable for
all $h\in H$ and every $\theta \in [0,1]$, (c.f.\ section~2.6 of \cite{12} for
more details).  Let $\DD_2(H)$ denote the class of $H$-valued $u(w)$ such that
$E_1|u|\subH\suptwo < \infty$.  The class of adapted square integrable random
variables is a closed subspace of $\DD_2(H)$ and will be denoted by
$\DD_2^a(H)$.  We will denote by $\hat u$ the projection of 
$u \in \DD_2(H)$ on $\DD_2^a (H)$, i.e.
\begin{equation}
\label{7.2}
E_1 |u-\hat u|\subH\suptwo = \inf_{v \in \DD_2^a(H)} E_1 |u-v|\subH\suptwo
\end{equation}
(this corresponds to the dual predictable projection in martingale theory).
Since $x$ is independent of $w$, and since
$\Hatx = E_1 (x|\sigma(y))$, then
\begin{equation}
\label{new7.3}
\hat x = (\Hatx)^{\widehat{}}
\end{equation}
i.e.\ if $x$ is not measurable on the $\sigma$-field induced by $y$, project,
first, $x$ on the $\sigma$-field generated by $y$ and then project on 
$\DD_2^a (H)$,
which is the same as replacing $u$ with $x$ in \eqref{7.2}.  Then (cf.\ e.g.
\cite{12})
\begin{equation}
\label{7.2plus1}
\ell(y) = \exp \left(\rho \delta \hat x - \frac{\rho^2}{2} |\hat x |\subH\suptwo\right)\,.
\end{equation}
By the same arguments as in \cite{1} or \cite{beta} and by the assumptions of
proposition~\ref{prop-51} 
\begin{equation}
\label{7.3}
I(X;Y) = \frac{\rho^2}{2} E |x-\hat{x}|\subH\suptwo
\,.
\end{equation}
\textbf{Remarks:}
(a)~The left hand side of \eqref{7.3} is independent of the choice of
$\pi_\bfcdot$ while $\widehat{x}$ does depend on the choice of 
$\pi_\bfcdot$.  Consequently, by \eqref{7.2plus1},
$E|x-\hat{x}|\subH\suptwo$ is independent of the particular choice of
$\pi_\bfcdot$.
(b)~The validity of \eqref{7.2plus1} and \eqref{7.3} is not restricted to the
case where $x$ is independent of $w$ (cf.\ \cite{beta} and \cite{12}).

By equations \eqref{5.2} and \eqref{7.3}, 
\begin{align*}
E_1 \log \ell (w) &=  \frac{\rho^2}{2} \: E\left( |x|\subH\suptwo - \half
|x-\hat{x}|\subH\suptwo\right)\\
& = \frac{\rho^2}{2} \: E \Bigl( - |\hat{x} |\subH\suptwo + 2(x, \hat{x})\subH \Bigr) 
\end{align*}
Hence
\begin{equation}
\label{7.4}
E_1 \log \ell(w) = \frac{\rho^2}{2} E|\hat{x}|\subH\suptwo
\,.
\end{equation}
We conclude the paper with the following remark.  Obviously
\begin{equation}
\label{new7.6}
E | \hat{x} |\subH\suptwo \le E|\Hatx|\subH\suptwo
\,.
\end{equation}
Hence by \eqref{7.4} and \eqref{6.3},
$$
E_1 \log \ell (w) \le \half E_1 | \nabla \log \ell (w)|\subH\suptwo
$$
or
\begin{equation}
\label{new7.7}
E_0  \ell (w)\log \ell (w) \le \half E_0\ell (w) | \nabla \log \ell (w)|\subH\suptwo
\,.
\end{equation}
Setting $f^2(w) = c \ell(w)$, $ c>0$ then
\begin{equation}
\label{7.6}
E_0 f^2 (w) 
\log |f(w)| \le E_0f^2 \cdot \log E_0^{\half} f^2 + E_0  |\nabla \log f|\subH\suptwo
\end{equation}
which is the logarithmic Sobolev inequality of L. Gross on Wiener space
(cf.\ e.g.\ section~9.2 of \cite{17} and the references therein).
Note, however, that \eqref{7.6}
is not the complete logarithmic Sobolev inequality since as derived above, it
holds only the for the case where $\ell(w)$ is the likelihood ratio associated
with $x+w$ where $x$ and $w$ are independent (and not for any nonnegative $\ell(w)$ for which
$E \ell(w) = 1$ cf. \cite{18}, \cite{19}).

Inequality \eqref{7.6} follows from the obvious inequality \eqref{new7.6} and
the \textit{equalities\/} derived earlier in this paper.  The question arises
whether a similar argument can yield \eqref{7.6} without the restriction on
$\ell(\om)$ to be generated by a signal plus \textit{independent\/} white
noise.  This seems to be a delicate problem;  the left hand side of
\eqref{new7.7} can be shown to be equal to the left hand side of
\eqref{new7.6} without the restriction that the signal and noise are
independent (\cite{beta}).  However it is not clear if the right hand side of
\eqref{new7.7} and \eqref{new7.6} are equal.

\end{document}